\documentclass[aos,MSNbibl,nameyear,dvips]{arximspdf}
\usepackage{graphicx}
\doi{10.1214/13-AOS1176} %kopijuoti is PTS
\volume{41}
\issue{6}
\pubyear{2013}
\firstpage{2979}
\lastpage{2993}

\makeatletter
\newtheorem{teo}{Theorem}
\newtheorem{lemma}{Lemma}%[section]

\newtheorem{cor}[teo]{Corollary}
\newproclaim{rem}{Remark}
\makeatother

\begin{document}
\begin{frontmatter}

\title{Note on distribution free testing for discrete~distributions}
\runtitle{Distribution free tests for discrete distributions}

\begin{aug}
\author[A]{\fnms{Estate} \snm{Khmaladze}\corref{}\ead[label=e1]{Estate.Khmaladze@vuw.ac.nz}}
\runauthor{E. Khmaladze}
\affiliation{Victoria University of Wellington}
\address[A]{School of Mathematics, Statistics and OR\\
Victoria University of Wellington\\
PO Box 600, Wellington\\
New Zealand\\
\printead{e1}} %adresu isvedimo komanda gale!
\end{aug}

% HISTORY:
\received{\smonth{8} \syear{2012}}
\revised{\smonth{9} \syear{2013}}

% ABSTRACT
%
\begin{abstract}
The paper proposes one-to-one transformation of the vector of
components $\{Y_{in}\}_{i=1}^m$ of Pearson's chi-square statistic,
\[
Y_{in}=\frac{\nu_{in} - np_i}{\sqrt{np_i}},\qquad i=1,\ldots, m,
\]
into another vector $\{Z_{in}\}_{i=1}^m$, which, therefore, contains
the same ``statistical information,'' but is asymptotically
distribution free. Hence any functional/test statistic based on
$\{Z_{in}\}_{i=1}^m$ is also asymptotically distribution free. Natural
examples of such test statistics are traditional goodness-of-fit
statistics from partial sums $\sum_{I\leq k} Z_{in}$.

The supplement shows how the approach works in the problem of
independent interest: the goodness-of-fit testing of power-law
distribution with the Zipf law and the Karlin--Rouault law as
particular alternatives.
\end{abstract}

% KEYWORDS
\begin{keyword}[class=AMS]
\kwd[Primary ]{62D05}
\kwd{62E20}
\kwd[; secondary ]{62E05}
\kwd{62F10}
\end{keyword}
\begin{keyword}
\kwd{Components of chi-square statistics}
\kwd{unitary transformations}
\kwd{parametric families of distributions}
\kwd{projections}
\kwd{power-law distributions}
\kwd{Zipf's law}
\end{keyword}

\end{frontmatter}

%s1 #&#
\section{Introduction}\label{sec1}
The main driver for this work was the need for a class of
distribution-free tests for discrete distributions. The basic step,
reported in Section~\ref{sec2} below, could have been made long ago,
maybe even soon after the publication of the classical papers of
\citet{Pear} and Fisher (\citeyear{Fish}, \citeyear{Fish2}). However, the tradition of using
the chi-square goodness-of-fit statistic became so widely spread, and
the point of view that, for discrete distributions, other statistics
``have to'' have their asymptotic distributions dependent on the
individual probabilities, became so predominant and ``evident,'' that
it required a strong impulse to examine the situation again. It came,
in this case, in the form of a question from Professor Ritei Shibata,
``Why is the theory of distribution-free tests for discrete
distributions so much more narrow than for continuous distributions?''
If it is true that sometimes a question is half of the answer, then
this is one such case.

We recall that for continuous distributions, the idea of the time
transformation $t=F(x)$ of \citet{Kol}, along with subsequent
papers of \citet{Sil} and \citet{WalWol}, was always
associated with a \textit{class} of goodness-of-fit statistics. The
choice of statistics invariant under this time transformation, at least
since the paper of \citet{AndDar}, became an accepted
\mbox{principle} in goodness-of-fit theory for continuous distributions. For
discrete distributions, however, everything is locked on a single
statistic, the chi-square goodness-of-fit statistic. It certainly is
true that in cases like the maximum likelihood statistic for
multinomial distributions [see, e.g., \citet{KenStu}] or
like the empirical likelihood [see, e.g., \citet{EinMcK}
and \citet{Owen}], the chi-square statistic appears as a natural
asymptotic object. Yet most of the time the choice of this statistic
comes as a deliberate choice of one particular asymptotically
distribution-free statistic. The idea of a \textit{class} of
asymptotically distribution free tests, to the best of our knowledge,
was never considered in any serious and systematic way.

This is a pity, because unlike the transformation $t=F(x)$, which is
basically a tool for one-dimensional time $x$, if we do not digress
onto the transformation of \citet{Ros} or spatial martingales of
\citet{Khm93}, the idea behind Pearson's chi-square test is
applicable to any measurable space. The potential of its generalization
seems, therefore, worth investigation.

We will undertake one such investigation in this paper. Namely, we will
obtain a transformation of the vector $Y_n$ of components of Pearson's
chi-square statistic (see below) into a vector $Z_n$, which will be
shown to be asymptotically distribution free. Therefore, any functional
based on $Z_n$ can be used as a statistic of an asymptotically
distribution-free test for the corresponding discrete distribution.
Thus the paper demonstrates, we hope, that the geometric insight behind
the papers of \citet{Pear} or \citet{Fish2} goes considerably further
than one goodness-of-fit statistic.

In the remaining part of this \hyperref[sec1]{Introduction} we present a typical result
of this paper. General results and other, may be more convenient, forms
of the transformation are given in the appropriate sections later on.

Let $p_1,\ldots,p_m$ be a discrete probability distribution; all $p_i >
0$ and $\sum_{i=1}^m p_i=1$. Denote $\nu_{1n}, \ldots, \nu_{mn}$ the
corresponding frequencies in a sample of size $n$, and consider the
vector $Y_n$ of components of the chi-square statistic
\[
Y_{in}=\frac{\nu_{in} - np_i}{\sqrt{np_i}},\qquad i=1,\ldots, m.
\]
Let $X=(X_1,\ldots,X_m)^T$ denote a vector of $m$ independent $N(0,1)$
random variables. As $n\to\infty$ the vector $Y_n$ has a limit
distribution of the zero-mean Gaussian vector $Y=(Y_1,\ldots,Y_m)^T$
such that
%
%e1 #&#
\begin{equation}
\label{chi} Y=X-\langle X,\sqrt p\rangle\sqrt p,
\end{equation}
where $\sqrt p$ denotes the vector $\sqrt p= (\sqrt p_1, \ldots, \sqrt
p_m)^T$. Here and below we use the notation $\langle a,b\rangle$ for
inner product of vectors $a$ and $b$ in ${\mathbb R}^m$:  $\langle
a,b\rangle=\sum_{i=1}^m a_i b_i $.

According to (\ref{chi}) the vector $Y$ is an orthogonal projection of
$X$ parallel to $\sqrt p$. Of course its distribution depends on $\sqrt
p$---it is only the sum of squares
\[
\langle Y,Y\rangle,
\]
which is chi-square distributed and hence has a distribution free from
$\sqrt p$. It is for this reason that we do not have any other
asymptotically distribution-free goodness-of-fit test for discrete
distributions except the chi-square statistic
\[
\langle Y_n,Y_n\rangle=\sum
_{i=1}^m \frac{(\nu_{in} -
np_i)^2}{np_i}.
\]
In particular, the asymptotic distribution of partial sums based on
$Y_{in}$, like
\[
\sum_{i=1}^k \frac{\nu_{in} - np_i} {\sqrt{np_i}} \quad\mbox{or}\quad\sum_{i=1}^k \frac{\nu_{in} - np_i}{\sqrt{n} },
\qquad k=1,2,\ldots, m,
\]
which would be discrete time analogues of the empirical process, will
certainly depend on $\sqrt p$, as will the asymptotic distribution of
statistics based on them. Here we would like to refer to paper of
\citet{Hen}, which advances the point of view that goodness-of-fit
tests for discrete distributions should be thought of as based on
empirical processes in discrete time, that is, on the partial sums on
the right. In the same vein, \citet{Chou} considered
quadratic functionals based on these partial sums, as direct analogues
of (weighted) omega-square statistics. We refer also to \citet{GolMor},
where tables for some quantiles of Kolmogorov--Smirnov
statistics from the partial sums are calculated in the parametric
problem, described in the supplementary material [\citet{Khm13}]. These papers illustrate the
dependence on the hypothetical distribution $p$ very clearly.

We do not know of many attempts to construct distribution-free tests
for discrete distributions, but one such, suggested in \citet{GreNik}, stands out for its simplicity and clarity: any discrete
distribution function $F_0$ can be replaced by a piece-wise linear
distribution function $\tilde F_0$ with the same values as $F_0$ at the
(nowhere dense) jump points of the latter; this opens up the
possibility to use time transformation $t=\tilde F_0(x)$ and thus
obtain distribution-free tests. However, without inquiring about the
consequences of implied additional randomization between the jump
points, this approach remains a one-dimensional tool.

In this paper we introduce a vector $Z_n=\{Z_{in}\}_{i=1}^m$ as
follows: let $r$ be the unit length ``diagonal'' vector with all
coordinates $1/\sqrt m$, and put
%
%e2 #&#
\begin{equation}
\label{z} Z_{n} = Y_{n} - \langle Y_n, r
\rangle\frac{1}{1- \langle\sqrt{p},r\rangle}(r - \sqrt p ).
\end{equation}
More explicitly,
\[
Z_{in} = \frac{\nu_{in} - np_i}{\sqrt{np_i}} - \frac{1}{\sqrt m} \sum
_{j=1}^m \frac{\nu_{jn} - np_j}{\sqrt{np_j}} \frac{1}{1 -
\sum_{j=1}^m \sqrt{p_j/m}} \biggl(
\frac{1}{\sqrt m} - \sqrt p_i \biggr).
\]
We will see that the following statement for $Z_n$ is true:

%==========================\mathbf Propsition====================
%
\begin{prop*}\label{notag}
Let $\mathbb{I}=(1,\ldots,1)^T$ denote the vector with all $m$
coordinates equal to $1$. The asymptotic distribution of $Z_n$ is that
of another, standard orthogonal projection
\[
Z\stackrel{d} {=}X-\langle X,r\rangle r = X- \frac{1}{m} \langle X,
\mathbb{I} \rangle\mathbb{I}
\]
and therefore any statistic based on $Z_n$ is asymptotically
distribution free. The transformation of $Y_n$ to $Z_n$ is one-to-one.
\end{prop*}

Thus the problem of testing $p$ is translated into the problem of
testing uniform discrete distribution of the same dimension $m$.

In particular, partial sums
\[
\sum_{i=1}^k Z_{in}, \qquad k
=1,2,\ldots, m,
\]
will asymptotically behave as a discrete time analog of the standard
Brownian bridge. On the other hand, since the transformation from $Y_n$
to $Z_n$ is one-to-one, $Z_n$ carries the same amount of statistical
information as $Y_n$.

For the proof of the proposition, see Theorem \ref{teoteo1} below. We
will see that this is not an isolated result, but one of several
possible results, and it follows from one particular point of view,
which is explained in the next section. We carry it on to the
parametric case in Section~\ref{sec3}.

%s2 #&#
\section{Pertinent unitary transformation}\label{sec2}
The idea behind the transformation (\ref{z}) can be explained as
follows: the problem with the vector $Y$ is that it projects a~standard
vector $X$ parallel to a specific vector, the vector $\sqrt p$. This
vector changes and with it changes the distribution of $Y$. However,
using an appropriate unitary operator, which incorporates $\sqrt p$,
one can ``turn'' $Y$ so that the result will be an orthogonal
projection parallel to a standard vector. One such standard vector can
be the vector $(1/\sqrt{m})\mathbb{I}$ above.

Slightly more generally, let $q$ and $r$ be two vectors of unit length
in $m$-dimensional space ${\mathbb R}^m$. Apart from obvious particular
choice of $r=(1/\sqrt m)\mathbb{I}$ and $q=\sqrt p= (\sqrt
p_1,\ldots,\sqrt p_m)^T$, we will consider other choices later on as
well. Denote by ${\mathcal L}={\mathcal L} (q,r)$ the $2$-dimensional
subspace of ${\mathbb R}^m$, generated by the vectors $q$~and~$r$, and
by ${\mathcal L^*} $ its orthogonal complement in ${\mathbb R}^m$. In
the lemma below we write $q_{\perp r}$ for the part of $q$ orthogonal
to $r$, and $r_{\perp q}$ for the part of $r$ orthogonal to $q$:
\[
q_{\perp r} = q - \langle q,r\rangle r,\qquad r_{\perp q} = r -
\langle q,r\rangle q
\]
and let $\mu=\|q_{\perp r}\|= \|r_{\perp q}\|$. Obviously, vectors $r$
and $q_{\perp r}/\mu$ form an orthonormal basis of ${\mathcal L}$ and
vectors $q$ and $r_{\perp q}/\mu$ form another orthonormal basis.
Consider
\[
U=r c^T + q_{\perp r} d^T/\mu
\]
with some
$c,d\in\mathcal{L}$, as a linear operator in ${\mathcal L}$.

%=================Lemma 1==============
%
%le1 #&#
\begin{lemma}\label{lem1}
\textup{(i)} The operator $U$ is unitary if and only
if the vectors $c$ and $d$ are orthonormal,
\[
\|c\|=\|d\|=1, \qquad\langle c,d\rangle=0.
\]

\textup{(ii)} The unitary operator $U$ maps $q$ to $r$,
\[
Uq=r,
\]
if and only if $c=q$ and $d=\pm r_{\perp q}/\mu$.

Altogether
\[
U= r q^T\pm\frac{1}{\mu^2} q_{\perp r} r_{\perp q}^T
\]
is the unitary operator in $\mathcal L$, which maps vector $q$ to
vector $r$. It also maps vector $r_{\perp q}$ to vector $\pm q_{\perp
r}$.
\end{lemma}

\begin{rem*}
In what follows in this section we will choose the sign $+$.

It is clear that if vector $x$ is orthogonal to $q$ and $r$, then
$Ux=0$. In other words, $U$ annihilates ${\mathcal L}^*$. Denote
$I_{{\mathcal L}^*}$ the projection operator parallel to $\mathcal L$,
so that it is the identity operator on ${\mathcal L}^*$ and annihilates
the subspace ${\mathcal L}$. Then the operator $I_{{\mathcal L}^*} + U$
is a unitary operator on ${\mathbb R}^m$. We use it to obtain our first
result.

Suppose vector $Y$ is projection of $X$, parallel to the vector $q$,
\[
Y=X-\langle X,q\rangle q.
\]
\end{rem*}

%
%==========================Theorem 1============
%
%th1 #&#
\begin{teo}\label{teoteo1}
\textup{(i)} The vector
%
%e3 #&#
\begin{equation}
\label{total} X' = (I_{{\mathcal L}^*} + U)X= X - \langle X, q
\rangle(q-r) - \langle X, r_{\perp q} \rangle\frac{1}{1- \langle
q,r\rangle}(r - q)
\end{equation}
is also a vector with independent $N(0,1)$ coordinates.

\textup{(ii)} The vector
%
%e4 #&#
\begin{equation}
\label{dmain} Z= (I_{{\mathcal L}^*} + U)Y = Y - \langle Y, r \rangle
\frac{1}{1-
\langle q,r\rangle}( r - q )
\end{equation}
is projection of $X'$ parallel to $r$,
\[
Z = X' - \bigl\langle X', r\bigr\rangle r.
\]
\end{teo}

\begin{pf}
(i)~By its definition, vector $Y$ is the orthogonal projection of $X$,
parallel to $q$. Therefore, if we project it further as
\[
R=Y-\langle Y,r_{\perp q}\rangle\frac{1}{\mu^2} r_{\perp q} = X -
\langle X, q\rangle q - \langle X,r_{\perp q}\rangle\frac{1}{\mu^2}
r_{\perp q},
\]
we will obtain the vector $R$ orthogonal to both $q$ and $r$, that is,
a vector in $\mathcal L^*$. If we apply operator $I_{{\mathcal L}^*}$
to $R$ it will not change, while $U$ will annihilate it, and thus
\begin{eqnarray*}
(I_{{\mathcal L}^*} + U)X&= & R + U \biggl(\langle X, q\rangle q +
\langle
X,r_{\perp q}\rangle\frac{1}{\mu^2} r_{\perp q} \biggr)
\\
&= & R + \langle X, q\rangle r + \langle X,r_{\perp q}\rangle
\frac{1}{\mu^2} q_{\perp
r}
\\
&=& X - \langle X, q \rangle(q-r) - \langle X, r_{\perp q} \rangle
\frac{1}{\mu^2} (r_{\perp q} - q_{\perp r} ).
\end{eqnarray*}
Noting that
\[
r_{\perp q} - q_{\perp r} = (r - q) \bigl(1+ \langle q,r\rangle\bigr)\quad\mbox{and}\quad\mu^2=1-\langle q,r\rangle^2,
\]
we obtain the right-hand side of (\ref{total}). Coordinates of $X'$ are
independent $N(0,1)$ random variables if the covariance matrix $EX'
X'^{T}$ is the identity matrix on ${\mathbb R}^m$. We have
\begin{eqnarray*}
EX' X'^{T} &= &(I_{{\mathcal L}^*} + U)EX
X^T (I_{{\mathcal L}^*} + U)^T = (I_{{\mathcal L}^*} + U)
\bigl(I_{{\mathcal L}^*} + U^T\bigr)
\\
&= & I_{{\mathcal
L}^*} + U U^T = I_{{\mathcal L}^*} + r r^T
+ \frac{1}{\mu^2} q_{\perp
r}q_{\perp r}^T = I.
\end{eqnarray*}

(ii) Note that the orthogonality property of $Y$, $\langle Y,
q\rangle=0$, implies that $\langle X,r_{\perp q}\rangle= \langle Y,
r\rangle$, and re-write (\ref{total}) as
\begin{eqnarray*}
X' &=& (I_{{\mathcal L}^*} + U)X = Y - \langle Y, r \rangle
\frac{1}{1-
\langle q,r\rangle}(r - q) + \langle X,q\rangle r.
\end{eqnarray*}
Also note that
\[
\bigl\langle X', r\bigr\rangle= \bigl\langle(I_{{\mathcal L}^*} +
U)X, r\bigr\rangle= \bigl\langle X, (I_{{\mathcal L}^*} + U)^Tr\bigr
\rangle= \langle X, q\rangle
\]
and so that $Z$ is indeed the projection of $X'$, we need
\[
Z = X' - \bigl\langle X', r\bigr\rangle r = Y- \langle
Y, r \rangle\frac
{1}{1- \langle q,r\rangle}(r - q).
\]\upqed
\end{pf}

The second statement of this theorem, together with the classical
statement $Y_n\stackrel{d}{\to} Y$, and the choice of
$r=(1,\ldots,1)/\sqrt m$ and $q=\sqrt p$, proves the proposition of the
\hyperref[sec1]{Introduction}.

The nature of the transformation and the proof given above does not
depend on a particular choice of the vector $r$ and is correct for any
$r$ of unit length. For example, we can choose $r=(1,0,\ldots,0)^T$.
Then the transformed vector $Z_n$ will have coordinates
%
%e5 #&#
\begin{equation}
\label{ind} Z_{in} = \frac{\nu_{in} - np_i}{\sqrt{np_i}} - \frac{\nu
_{1n} -
np_1}{\sqrt{np_1}}
\frac{1}{1 - \sqrt{p_1}} (\delta_{1i} - \sqrt{p_i} )
\end{equation}
or
\[
Z_{1n}=0,\qquad
Z_{in} = \frac{\nu_{in} - np_i}{\sqrt{np_i}} -
\frac{\nu_{1n} - np_1}{\sqrt{np_1}} \frac{1}{1 - \sqrt{p_1}} \sqrt
{p_i},
\qquad i=2,\ldots,m.
\]
As a corollary of the previous theorem we obtain a vector with very
simple asymptotic behavior.

%======================Cor 2, to (0,X,..., X) ===================
%
%co2 #&#
\begin{cor}\label{cor2} If $Y_n\stackrel{d}{\to} Y=X-\langle
X,\sqrt{p}\rangle\sqrt{p}$, then for the vector $Z_n$ defined in
(\ref{ind}) we have
\[
Z_n\stackrel{d} {\to} (0,X_2,\ldots,X_m)^T.
\]
\end{cor}

To find the asymptotic distribution of statistics based on this choice
of $Z_n$ may be more convenient than in the previous case. Yet the
relationship between the two is one-to-one.

It is often the case that the probabilities $p_1,\ldots,p_m$ depend on
a parameter, which has to be estimated from observed frequencies. This
case needs additional consideration which we defer to the next section.
However, there are also cases when the hypothetical probabilities are
fixed, or the value of the parameter is estimated from previous
samples, and therefore needs to be treated as a given. In these cases
Theorem \ref{teoteo1} is directly applicable.

One important case of this type is the two-sample problem. Namely, let
events, labeled by $i=1,2,\ldots,m$, be basically as above, and let
$\nu'_{1n'}, \ldots, \nu'_{mn'}$ and $\nu''_{1n''}, \ldots,
\nu''_{mn''}$ be frequencies of these events in two independent samples
of size $n'$~and~$n''$, respectively. Let $\mu_{1}, \ldots, \mu_{m}$
denote the frequencies in the pooled sample of size $n=n'+n''$. Then
the normalized differences
\[
Y'_{in}=\frac{\nu_{in'}' - n' \mu_{i}/n }{\sqrt{n' \mu
_{i}/n}},\qquad i=1,\ldots,m,
\]
are the components of the two sample chi-square statistic: the sum of
their squares is the statistic. Conditions which guarantee convergence
of the vector $Y'_n$ of these differences in distribution to the vector
$Y$ are well known; see, for example, \citet{Rao}, or \citet{EinKhm}
and references therein. Then it follows from Theorem
\ref{teoteo1} that under these conditions the vector $Z'_n$ with
coordinates
\begin{eqnarray*}
Z'_{in} &=& \frac{\nu'_{in'}
- n' \mu_{i}/n }{\sqrt{n' \mu
_{i}/n}}
\\
&&{} - \frac{1}{\sqrt m} \sum _{j=1}^m \frac{\nu'_{jn} - n'
\mu_{j}/n }{\sqrt{n' \mu_{j}/n }}
\frac{1}{1 + \sum_{j=1}^m \sqrt{\mu_{j}/nm}} \biggl( \frac
{1}{\sqrt m} + \sqrt{\frac{\mu_{i}}{n}}
\biggr)
\end{eqnarray*}
converges in distribution to vector $X-\langle X,\mathbb{I}\rangle
\mathbb{I}/m$ and, hence, is asymptotically distribution free. To show
this result one needs only to choose as $q$ the vector $(\sqrt
{\mu_{1}/n}, \ldots, \sqrt{\mu_{m}/n})^T$ in Theorem \ref{teoteo1}
above. Corollary \ref{cor2} suggests another choice of the transformed
vector with coordinates
\[
Z_{in} = \frac{\nu'_{in'} - n' \mu_{i}/n }{\sqrt{n' \mu_{i}/n}} -
\frac{\nu'_{1n} - n' \mu_{1}/n }{\sqrt{n' \mu_{1}/n }} \frac{1}{1 +
\sqrt{\mu_{1}/n}}
\sqrt{\frac{\mu_{i}}{n}},\qquad i=2,\ldots,m
\]
with also simple asymptotic behavior.

%========================= Sec. 3 parametric discrete ===============

%s3 #&#
\section{The case of estimated parameters}\label{sec3}
%============================discrete parametric
We will now see that the pivotal property of $Y_n$ to behave as
asymptotically orthogonal projection of $X$ remains true for components
of chi-square statistic with estimated parameter.

Indeed, if the hypothetical probabilities depend on a
$\kappa$-dimensional parameter, $p_i=p_i(\theta)$, which is estimated
via maximum likelihood or minimum chi-square, then the statistic
\[
\sum_{i=1}^m \frac{(\nu_{in} - np_i(\hat\theta_n))^2}{np_i(\hat
\theta_n)}
\]
has chi-square distribution with $m-1-k$ degrees of freedom; see
extensive review of this matter in \citet{Stig}, Chapter~19.
Notwithstanding great convenience of this result, note, however, that
the asymptotic distribution of the vector $\hat Y_n$ itself, with
%
%e6 #&#
\begin{equation}
\label{3-1} \hat Y_{in} = \frac{\nu_{in} - np_i(\hat\theta_n)}{\sqrt
{np_i(\hat
\theta_n)}},
\end{equation}
depends, under hypothesis, not only on the probabilities $p_i(\theta)$
at the true value of~$\theta$, but also on their derivatives in
$\theta$. Therefore, the limit distribution of statistics from $\hat
Y_n$ in general will depend on the hypothetical parametric family and
on the value of the parameter.

At the same time, it is well known since long ago [see, e.g.,
\citet{Cra}, Chapter~20; a modern treatment can be found in
\citet{vdV}] that under mild assumptions the maximum likelihood
(and minimum chi-square) estimator possesses asymptotic expansion of
the form
\begin{eqnarray*}
{\sqrt n} (\hat\theta_n - \theta) &=& \Gamma^{-1} \sum
_{i=1}^m Y_{in}
\frac{\dot p_i(\theta)}{\sqrt{p_i(\theta)}} + o_P(1),
\end{eqnarray*}
where $\dot p_i(\theta)$ denotes the $\kappa$-dimensional vector of
derivatives of $p_i(\theta)$ in $\theta$ and
\[
\Gamma=\sum_{i=1}^m \frac{\dot p_i(\theta) \dot p_i(\theta)^T}
{p_i(\theta)}
\]
denotes the
$\kappa\times\kappa$ Fisher information matrix. At the same time, the
expansion
\[
\label{yhat1} \hat Y_{in} = Y_{in} - \frac{\dot p_i(\theta)^T}{\sqrt
{p_i(\theta
)}} {
\sqrt n} (\hat\theta_n - \theta) + o_P(1)
\]
is also true. Combining these two expansions, one obtains
%
%e7 #&#
\begin{equation}
\label{yhat} \hat Y_{in} = Y_{in} - \frac{\dot p_i(\theta)^T}{\sqrt
{p_i(\theta
)}}
\Gamma^{-1} \sum_{i=1}^m
Y_{in} \frac{\dot p_i(\theta
)}{\sqrt{p_i(\theta)}} + o_P(1).
\end{equation}
Use the notation
\[
\hat q_i= \Gamma^{-1/2}\frac{\dot p_i(\theta)}{\sqrt{p_i(\theta)}},\qquad
i=1,\ldots,m
\]
and remember that
\[
\sum_{i=1}^m\sqrt{p_i(
\theta)} \frac{\dot p_i(\theta)^T}{\sqrt{p_i(\theta)}} = 0,
\]
that is, that the vectors in $i$, which form $\dot p/{\sqrt p}$, are
orthogonal to the vector ${\sqrt p}$. Therefore all $\kappa$
coordinates of $q_i$ form, in $i$, vectors which are orthonormal and
orthogonal to the vector $\sqrt{p(\theta)}$. Together with (\ref{chi})
this implies the convergence in distribution of $\hat Y_{n}$ to
Gaussian vector
%
%e8 #&#
\begin{equation}
\label{limyhat} \hat Y = X - \langle X,\sqrt p\rangle{\sqrt p} -
\langle X, \hat
q\rangle{\hat q}.
\end{equation}

It is easily seen that expression (\ref{limyhat}) describes $\hat Y$ as
an orthogonal projection of~$X$ parallel to vectors ${\sqrt p}$ and
$\dot p/{\sqrt p}$; see \citet{Khm} for an analogous description
of empirical processes. Using this description, we can extend the
method of Section~\ref{sec2} to the present situation.

Indeed, let us assume from now on that $\kappa=1$, which will make the
presentation more transparent. Having two vectors, $q=\sqrt{p(\theta)}$
and $\hat q$, which determine the asymptotics of $\hat Y_n$, let us
choose now a standard vector $r$ of unit length and another vector,
$\hat r$, also of unit length and orthogonal to $r$. Heuristically, one
may think of it as a normalized ``score function'' for some
``standard'' family around $r$. For example, choose
$r=(1/\sqrt{m})\mathbb I$ and choose any unit vector, such that
$\sum_{i=1}^m \hat r_i =0$. Two such choices, we think, will be
particularly useful: for $m$ even,
\[
\frac{1}{\sqrt m}(1,\ldots,1, - 1, \ldots, - 1)^T
\]
or
\[
\frac{1}{\sqrt m}(1,\ldots,1, - 1, \ldots, - 1, 1, \ldots, 1)^T
\]
with the ``plateau'' of $-1$s taken $m/2$-long, and for $m$ odd put,
say, the last coordinate equal $0$.

Whatever the choice of $\hat r$, suppose we chose and fixed it. It is
obvious that the vector
%
%e9 #&#
\begin{equation}
\label{zhat} \hat Z=X - \langle X, r\rangle r - \langle X, \hat r\rangle
\hat r
\end{equation}
has a distribution totally unconnected, and hence free from the
parametric family~$p(\theta)$. Consider now the subspace $\hat{\mathcal
L}={\mathcal L} (q, \hat q, r, \hat r)$. We do not need to insist that
it is a 4-dimensional subspace, but typically it is, at least, as far
as we have freedom in $\hat r$. Let $\hat{\mathcal L}{}^*$ denote the
orthogonal complement of $\hat{\mathcal L}$ to ${\mathbb R}^m$. Two
bases of the space $\hat{\mathcal L}$ will be useful: one is formed by
$r,\hat r, b_3, b_4$ where $b_3$ and $b_4$ are re-arrangements of $q$
and $\hat q$, which are orthonormal and orthogonal to $r$~and~$\hat r$;
the other is formed by $q,\hat q, a_3, a_4$ where $a_3$~and~$a_4$ are,
re-arrangements of $r$~and~$\hat r$, which are orthonormal and
orthogonal to $q$ and $\hat q$. We will consider particular forms of
these vectors later on.

%==========================Lemma 2==============
%
%le2 #&#
\begin{lemma}\label{lem2}
The operator
\[
\hat U=r q^T + \hat r \hat q{}^T + b_3
a_3^T + b_4a_4^T
\]
is a unitary operator on $\hat{\mathcal L}$ and such that
\[
\hat U q=r, \qquad\hat U \hat q = \hat r.
\]
\end{lemma}

%===========================Theorem 3==why 3?==============
%
%th3 #&#
\begin{teo}\label{therm3}
Under convergence in distribution of the vector $\hat Y_n$ with
coordinates (\ref{3-1}) to the Gaussian vector $\hat Y$ given by
(\ref{limyhat}), the vector
%
%e10 #&#
\begin{equation}
\label{teo3} \hat Z_n = \hat Y_n - \langle\hat
Y_n, a_3\rangle(a_3 - b_3) -
\langle\hat Y_n, a_4\rangle(a_4 -
b_4)
\end{equation}
converges in distribution to the Gaussian vector $\hat Z$ given by
(\ref{zhat}). Therefore, any statistic based on $Z_n$ is asymptotically
distribution free.
\end{teo}

\begin{pf}
Let $\mathcal{ \hat L }{}^*$ be orthogonal complement of the subspace
$\hat{\mathcal{L}}$ in ${\mathbb R}^m$ and let $\hat I$ be projector on
the $\hat{\mathcal{L}}{}^*$. We need to verify two things: (a)~that the
vector $\hat Z$ can be obtained as
\[
\hat Z = (\hat I +\hat U)\hat Y
\]
and (b) that its explicit form is as given in the theorem. We show (a)
slightly differently from what was done in Theorem \ref{teoteo1}.
Namely, recall that the covariance operator of $\hat Y$ is the
projector $E\hat Y \hat Y{}^T = I - q q^T - \hat q \hat q{}^T$, where $I$
stands for an identity operator on ${\mathbb R}^m$, and consider the
covariance operator of $(\hat I +\hat U)\hat Y$:
\[
E (\hat I +\hat U)\hat Y \hat Y{}^T (\hat I +\hat U)^T = (
\hat I +\hat U) \bigl(I - q q^T - \hat q \hat q{}^T\bigr) (
\hat I +\hat U)^T.
\]
However, $(\hat I +\hat U)I (\hat I + \hat U)^T =I$ while $(\hat I
+\hat U)q=r$ and $(\hat I +\hat U)\hat q=\hat r$. This implies that
\[
(\hat I +\hat U) \bigl(I - q q^T - \hat q \hat q{}^T\bigr)
(\hat I +\hat U)^T = I-r r^T - \hat r \hat
r{}^T,
\]
which is the covariance operator of $\hat Z$.

To show (b) use the basis $q,\hat q, a_3, a_4$ and the orthogonality of
$\hat Y$ to $q$ and $\hat q$ to find that the projection of $\hat Y$
on $\mathcal{\hat L}$ can be written as
\[
\langle\hat Y, a_3\rangle a_3 + \langle\hat Y,
a_4\rangle a_4
\]
and therefore the difference $\hat Y - \langle\hat Y, a_3\rangle a_3 -
\langle\hat Y, a_4\rangle a_4 $ will remain unchanged by the operator
$\hat I$. At the same time $\hat U a_3=b_3$ and $\hat U a_4=b_4$. This
leads to the following form of our transformed vector $\hat Z$:
\[
(\hat I +\hat U)\hat Y = \hat Y - \langle\hat Y, a_3
\rangle(a_3 - b_3) - \langle\hat Y, a_4
\rangle(a_4 - b_4).
\]\upqed
\end{pf}

With regard to practical applications, there are several natural
choices of vectors $a_3, a_4$. For example, denote $r_{\perp q \hat q}$
the part of $r$ orthogonal to both $q$ and $\hat q$, and choose
\[
a_3 = \frac{1}{\|r_{\perp q \hat q}\|} r_{\perp q \hat q}= \frac{1}{\|
r_{\perp q \hat q}\|}
\bigl(r-\langle r,q\rangle q - \langle r,\hat q\rangle\hat q \bigr)
\]
and, similarly, choose $a_4$ as
\[
a_4 = \frac{1}{\|\hat r_{\perp r q \hat q}\|} \hat r_{\perp r q \hat
q}= \frac{1}{\|r_{\perp r q \hat q}\|}
\bigl(\hat r-\langle\hat r,q\rangle q - \langle\hat r,\hat q\rangle\hat
q -
\langle\hat r,a_3\rangle a_3\bigr).
\]
In dual way, we can choose specific $b_3$ and $b_4$ as
\[
b_3 = \frac{1}{\|q_{\perp r \hat r}\|} q_{\perp r \hat r}= \frac{1}{\|
q_{\perp r \hat r}\|}
\bigl(q -\langle q,r\rangle r - \langle q,\hat r\rangle\hat r \bigr)
\]
and
\[
b_4 = \frac{1}{\|\hat q_{\perp q r \hat r}\|} \hat q_{\perp q r \hat
r}= \frac{1}{\|q_{\perp q r \hat r}\|}
\bigl(\hat q-\langle\hat q, r\rangle r - \langle\hat q,\hat r\rangle
\hat r -
\langle\hat q,b_3\rangle b_3\bigr).
\]

A more symmetric choice would be
\[
a_3=\frac{1}{\sqrt2}\frac{1}{\sqrt{1+\rho}} \biggl(\frac{1}{\|r_{\perp
q \hat q}\|}
r_{\perp q \hat q} + \frac{1}{\|\hat r_{\perp q \hat q}\|
} \hat r_{\perp q \hat q}\biggr)
\]
and
\[
a_4=\frac{1}{\sqrt2} \frac{1}{\sqrt{1- \rho}} \biggl(\frac{1}{\|
r_{\perp q \hat q}\|}
r_{\perp q \hat q} - \frac{1}{\|\hat r_{\perp q
\hat q}\|} \hat r_{\perp q \hat q}\biggr),
\]
where $\rho$ is correlation\vspace*{-1pt} coefficient between $r_{\perp q \hat q}$
and $\hat r_{\perp q \hat q} $. Note that in both cases the inner
products $\langle\hat Y, a_3\rangle$ and $\langle\hat Y, a_4\rangle$
become linear combinations of just $\langle\hat Y, r\rangle$ and
$\langle\hat Y, \hat r\rangle$. For the last, symmetric choice, for
example, they are
\[
\frac{1}{\sqrt2}\biggl(\frac{1}{\|r_{\perp q \hat q}\|} \langle\hat Y,
r\rangle\pm
\frac{1}{\|\hat r_{\perp q \hat q}\|} \langle\hat Y, \hat r\rangle\biggr),
\]
respectively.

Although the choice of $r=(1/\sqrt m)\mathbb I$ is a natural one, the
different choice of the vectors $r$ and $\hat r$ leads to simpler form
of the transformed vector with convenient and simple asymptotic
distribution. Namely, let $r=(1,0,\ldots, 0)^T$ and $\hat
r=(0,1,0,\ldots,0)^T$. Then $\langle\hat Y, r\rangle$ and $\langle\hat
Y, \hat r\rangle$ become
\[
\frac{1}{\sqrt2} \frac{1}{\sqrt{1\pm\rho}} \biggl(\frac{1}{\sqrt{1 -
q_1^2 - \hat q_1^2 }} \hat
Y_1 \pm\frac{1}{ \sqrt{1 - q_2^2 - \hat
q_2^2 }} \hat Y_2 \biggr),
\]
respectively, with
\[
\rho=\frac{-q_1 q_2 - \hat q_1 \hat q_2}{\sqrt{1 - q_1^2 - \hat
q_1^2 } \sqrt{1 - q_2^2 - \hat q_2^2 }}.
\]
The form of vectors $a_3, a_4, b_3$ and $b_4$ also becomes simpler.
Similar to Corollary~\ref{cor2}, we have the following:

%===========================Cor 4, \to(0,0,X,...,X) ============
%
%co4 #&#
\begin{cor}\label{corh2}
If $r=(1,0,\ldots,0)^T$ and $\hat r=(0,1,0,\ldots,0)^T$ and if $\hat
Y_n \stackrel{d}{\to} \hat Y$ with $\hat Y$ described in
(\ref{limyhat}), then for the vector $\hat Z_n$ described in the
Theorem \ref{therm3}, we have
\[
\hat Z_n\stackrel{d} {\to} \hat Z=(0,0, X_3,\ldots,
X_m)^T,
\]
where $X_3,\ldots, X_m$ are independent and $N(0,1)$-distributed.
\end{cor}

\begin{rem*}
Although explicit coordinate representation through vectors $a_3$,
$a_4$, $b_3$, $b_4$ is useful in several ways, another representation
may be simpler, especially when more than one parameter is present. Let
us start with notation
\[
U_{q,r} = I - \frac{2}{\|r-q\|^2}(r-q) (r-q)^T.
\]
This is a unitary operator in ${\mathbb R}^m$, which maps $q$ into $r$
and $r$
into $q$, while any vector orthogonal to $r$ and $q$ is mapped into
itself. Note that $\|r-q\|$ is Hellinger distance between distributions
given by probabilities $(r_1^2,\dots, r_m^2)$ and $(q_1^2,\dots,
q_m^2)$ and that
\[
\|r-q\|^2 = 2\bigl(1-\langle q,r\rangle\bigr).
\]
We thus see that $U_{q,r}$ is simply a shorter notation for the
operator $I_{{\mathcal L}^*} + U$ of Section~\ref{sec2}. Now consider
an image $\tilde q= U_{q,r} \hat q$ of $\hat q$. This vector is
orthogonal to~$r$. Consider another operator $U_{\tilde q, \hat r}$.
Since both $\tilde q$ and $\hat r$ are orthogonal to $r$, this operator
will leave $r$ unchanged,\vspace*{-1pt} while mapping $\tilde q$ to $\hat r$. The
product $U_{\tilde q,\hat r} U_{q,r}$ will be another form of the
operator $\hat I + \hat U$, and (\ref{teo3}) can be written
as\looseness=-1
\[
\hat Z_n=U_{\tilde q,\hat r} U_{q,r} \hat Y_n.
\]\looseness=0
This recursive representation can obviously be extended for any $\kappa
>1$.
\end{rem*}

%============================ Sec 4 =====================
%===============================numerical illustrations===============
%s4 #&#
\section{On numerical illustrations}\label{sec4} One would hope that numerical
verification of the whole approach will be attempted in the future.
This will require a substantial amount of time and more room than the
present paper could allow. We also stress that this paper does not
advocate any particular test; its aim is to provide a~satisfactory
foundation on which various goodness-of-fit tests can be based.
However, in the supplementary material [\citet{Khm13}] we tried the approach on a
testing problem of independent interest: goodness-of-fit testing of the
power-law distributions with the Zipf law and the Karlin--Rouault law
as alternatives. We show some illustrations of how particular test
statistics based on partial sums of $Y_{in}$ and partial sums of
$Z_{in}$ perform in this problem.

In this section we restrict ourselves with one numerical illustration
of how quickly the asymptotic distribution freeness of vector $\hat
Z_n$ of (\ref{teo3}) start manifesting itself for finite $n$. For this
we considered three different choices of $p_1,\dots,p_m$ of the same
$m=10$. As the first choice we picked these probabilities at random: 9
uniform random variables have been generated once and the resulting
uniform spacings were used as these probabilities; as the second and
third choices we used increments $\Delta F(i/10)$, $i=1,\ldots,10$, of
beta distribution function with a bell shaped density, with parameters
3 and 3, and then with $J$-shaped density, with parameters 0.8 and 1.5.

%f1 #&#
\begin{figure}

\includegraphics{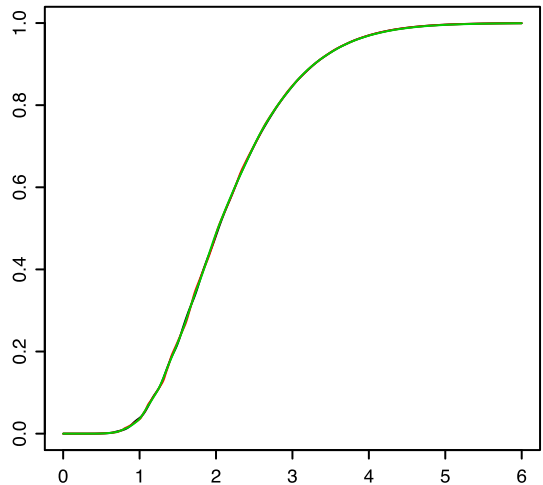}

\caption{Distribution functions of the statistic $d_{mn}^Z$ for three
different discrete distributions, as described in the text. 10,000
simulations of samples of size $n=200$ have been used. The dimension of
the discrete distributions (number of different events) was $m=10$.}\label{figdiscrete}
\end{figure}

From each of these distributions we generated 10,000 samples of size
$n=200$, and for each sample calculated a discrete version of the
Kolmogorov--Smirnov statistic
\[
d_{mn}^Z = \max_{1\leq k\leq m} \biggl|\sum
_{j\leq k} Z_{in} \biggr|.
\]
Figure~\ref{figdiscrete} shows three graphs of the resulting empirical
distribution functions.

In our choice of $n$ we tried to achieve what is typically required for
an application of Pearson's chi-square statistics, that all $np_i$ will
be at least $10$. Otherwise\vadjust{\goodbreak} we tried to choose $n$ not large. For
$n=200$ the requirement $np_i\geq10$ was not strictly satisfied, and in
the last two cases we had about three cells with $np_i$ about 5. This
could have somewhat spoiled the asymptotic result, but has not. If the
three graphs are not very distinct, that is because for all three cases
they are very close. Our statistic $d_{mn}^Z$ indeed looks distribution
free.

%s5 #&#
\section{Acknowledgment}\label{sec5}
For numerical results of the last section and
in the supplementary material [\citet{Khm13}] I am indebted to Boyd Anderson and Thuong Nguyen, and
also to Dr Ray Brownrigg.

\begin{supplement}\label{supp}
\stitle{Supplement: Distribution free Kolmogorov--Smirnov and
Cram\'er--von Mises tests for power-law distribution}
\slink[doi]{10.1214/13-AOS1176SUPP} \sdatatype{.pdf}
\sfilename{aos1176\_supp.pdf}
\sdescription{We compare asymptotic
behavior of the two classical goodness-of-fit tests based on partial
sums of $Y_{in}$'s and their distribution free transformations
$Z_{in}$'s and show their power under
Zipf's law and under Karlin--Rouault law as alternatives.} %For some
%reason the title and the description of the supplement appear at the
%end
\end{supplement}

%suskaldyti doi

% imsref loaded by linak, 2013-11-26 09:21:44

\printaddresses


\begin{thebibliography}{29}
% BibTex style file: ims.bst, 2013-01-28
% Default style options (sort=0,type=number).
% Used options (sort=1,type=nameyear).
%b1 #&#
\bibitem[\protect\citeauthoryear{Anderson and Darling}{1952}]{AndDar}
\begin{barticle}[mr]
\bauthor{\bsnm{Anderson},~\bfnm{T.~W.}\binits{T.~W.}} \AND
  \bauthor{\bsnm{Darling},~\bfnm{D.~A.}\binits{D.~A.}}
(\byear{1952}).
\btitle{Asymptotic theory of certain ``goodness of fit'' criteria based on
  stochastic processes}.
\bjournal{Ann. Math. Statistics}
\bvolume{23}
\bpages{193--212}.%
\bid{issn={0003-4851}, mr={0050238}}%
\bptok{imsref}%
\end{barticle}%
\endbibitem

%%b2 #&#
%(\byear{2002}).

%b3 #&#
\bibitem[\protect\citeauthoryear{Choulakian, Lockhart and Stephens}{1994}]{Chou}
\begin{barticle}[mr]
\bauthor{\bsnm{Choulakian},~\bfnm{V.}\binits{V.}},
  \bauthor{\bsnm{Lockhart},~\bfnm{R.~A.}\binits{R.~A.}} \AND
  \bauthor{\bsnm{Stephens},~\bfnm{M.~A.}\binits{M.~A.}}
(\byear{1994}).
\btitle{Cram\'er--von {M}ises statistics for discrete distributions}.
\bjournal{Canad. J. Statist.}
\bvolume{22}
\bpages{125--137}.
\bid{doi={10.2307/3315828}, issn={0319-5724}, mr={1271450}}
\bptnote{check year}%
\bptok{imsref}%
\end{barticle}
\endbibitem

%%b4 #&#
%  \bauthor{\bsnm{Shalizi},~\bfnm{Cosma~Rohilla}\binits{C.~R.}} \AND
%  \bauthor{\bsnm{Newman},~\bfnm{M.~E.~J.}\binits{M.~E.~J.}}
%(\byear{2009}).

%b5 #&#
\bibitem[\protect\citeauthoryear{Cram{\'e}r}{1946}]{Cra}
\begin{bbook}[mr]
\bauthor{\bsnm{Cram{\'e}r},~\bfnm{Harald}\binits{H.}}
(\byear{1946}).
\btitle{Mathematical {M}ethods of {S}tatistics}.
\bpublisher{Princeton Univ. Press}, \blocation{Princeton}.
\bid{mr={0016588}}
\bptok{imsref}%
\end{bbook}
\endbibitem\vadjust{\goodbreak}

%b6 #&#
\bibitem[\protect\citeauthoryear{Einmahl and Khmaladze}{2001}]{EinKhm}
\begin{bincollection}[mr]
\bauthor{\bsnm{Einmahl},~\bfnm{J.~H.~J.}\binits{J.~H.~J.}} \AND
  \bauthor{\bsnm{Khmaladze},~\bfnm{E.~V.}\binits{E.~V.}}
(\byear{2001}).
\btitle{The Two-Sample Problem in {$\Bbb R\sp m$} and Measure-Valued Martingales}.
%In \bbooktitle{State of the Art in Probability and Statistics ({L}eiden,  1999)}.
\bseries{Institute of Mathematical Statistics Lecture Notes---Monograph Series}
\bvolume{36}
\bpages{434--463}.
\bpublisher{IMS}, \blocation{Beachwood, OH}.
\bid{doi={10.1214/lnms/1215090082}, mr={1836574}}
\bptok{imsref}%
\end{bincollection}
\endbibitem

%b7 #&#
\bibitem[\protect\citeauthoryear{Einmahl and McKeague}{1999}]{EinMcK}
\begin{barticle}[mr]
\bauthor{\bsnm{Einmahl},~\bfnm{John H.~J.}\binits{J.~H.~J.}} \AND
  \bauthor{\bsnm{McKeague},~\bfnm{Ian~W.}\binits{I.~W.}}
(\byear{1999}).
\btitle{Confidence tubes for multiple quantile plots via empirical likelihood}.
\bjournal{Ann. Statist.}
\bvolume{27}
\bpages{1348--1367}.
\bid{doi={10.1214/aos/1017938929}, issn={0090-5364}, mr={1740108}}
\bptok{imsref}%
\end{barticle}
\endbibitem

%b8 #&#
\bibitem[\protect\citeauthoryear{Fisher}{1922}]{Fish}
\begin{barticle}[auto:STB|2013/10/14|10:36:11]
\bauthor{\bsnm{Fisher},~\bfnm{R.~A.}\binits{R.~A.}}
(\byear{1922}).
\btitle{On the interpretation of $\chi^2$ from contingency tables, and the
  calculation of~$P$}.
\bjournal{J. R. Stat. Soc.}
\bvolume{85}
\bpages{87--94}.
\bptok{imsref}%
\end{barticle}
\endbibitem

%b9 #&#
\bibitem[\protect\citeauthoryear{Fisher}{1924}]{Fish2}
\begin{barticle}[auto:STB|2013/10/14|10:36:11]
\bauthor{\bsnm{Fisher},~\bfnm{R.~A.}\binits{R.~A.}}
(\byear{1924}).
\btitle{Conditions under which $\chi^2$ measures the discrepancy between
  observation and hypothesis}.
\bjournal{J. R. Stat. Soc.}
\bvolume{87}
\bpages{442--450}.
\bptok{imsref}%
\end{barticle}
\endbibitem

%b10 #&#
\bibitem[\protect\citeauthoryear{Goldstein, Morris and Yen}{2004}]{GolMor}
\begin{barticle}[auto:STB|2013/10/14|10:36:11]
\bauthor{\bsnm{Goldstein},~\bfnm{M.~L.}\binits{M.~L.}},
  \bauthor{\bsnm{Morris},~\bfnm{S.~A.}\binits{S.~A.}} \AND
  \bauthor{\bsnm{Yen},~\bfnm{G.~G.}\binits{G.~G.}}
(\byear{2004}).
\btitle{Problems with fitting the power-law distributions}.
\bjournal{Eur. Phys. J. B}
\bvolume{41}
\bpages{255--258}.
\bptok{imsref}%
\end{barticle}
\endbibitem

%b11 #&#
\bibitem[\protect\citeauthoryear{Greenwood and Nikulin}{1996}]{GreNik}
\begin{bbook}[mr]
\bauthor{\bsnm{Greenwood},~\bfnm{Priscilla~E.}\binits{P.~E.}} \AND
  \bauthor{\bsnm{Nikulin},~\bfnm{Mikhail~S.}\binits{M.~S.}}
(\byear{1996}).
\btitle{A Guide to Chi-Squared Testing}.
\bpublisher{Wiley}, \blocation{New York}.
\bid{mr={1379800}}
\bptok{imsref}%
\end{bbook}
\endbibitem

%%b12 #&#
%  \bauthor{\bsnm{Chibisov},~\bfnm{D.~M.}\binits{D.~M.}}
%(\byear{1979}).
%In \bbooktitle{Contributions to Statistics}

%b13 #&#
\bibitem[\protect\citeauthoryear{Henze}{1996}]{Hen}
\begin{barticle}[mr]
\bauthor{\bsnm{Henze},~\bfnm{Norbert}\binits{N.}}
(\byear{1996}).
\btitle{Empirical-distribution-function goodness-of-fit tests for discrete models}.
\bjournal{Canad. J. Statist.}
\bvolume{24}
\bpages{81--93}.
\bid{doi={10.2307/3315691}, issn={0319-5724}, mr={1394742}}
\bptnote{check year}%
\bptok{imsref}%
\end{barticle}
\endbibitem

%b14 #&#
\bibitem[\protect\citeauthoryear{Kendal and Stuart}{1963}]{KenStu}
\begin{bmisc}[auto:STB|2013/10/14|10:36:11]
\bauthor{\bsnm{Kendal},~\bfnm{M.}\binits{M.}} \AND
  \bauthor{\bsnm{Stuart},~\bfnm{A.}\binits{A.}}
(\byear{1963}).
\bhowpublished{\textit{The advanced theory of statistics}, \textit{Vol.~2}. C.~Griffin, London. Re-printed as
  Vol.~2A.---Classical Inference and the Linear Model in 2009}.
\bptok{imsref}%
\end{bmisc}
\endbibitem

%b15 #&#
\bibitem[\protect\citeauthoryear{Khmaladze}{1979}]{Khm}
\begin{barticle}[auto:STB|2013/10/14|10:36:11]
\bauthor{\bsnm{Khmaladze},~\bfnm{E.~V.}\binits{E.~V.}}
(\byear{1979}).
\btitle{The use of omega-square tests for testing parametric hypotheses}.
\bjournal{Theory Probab. Appl.}
\bvolume{v.XXIV}
\bpages{283--302}.
\bptok{imsref}%
\end{barticle}
\endbibitem

%%b16 #&#
%(\byear{1989}).
%  Report, Amsterdam}.

%b17 #&#
\bibitem[\protect\citeauthoryear{Khmaladze}{2013}]{Khm13}
\begin{bmisc}[auto:STB|2013/10/14|10:36:11]
\bauthor{\bsnm{Khmaladze},~\bfnm{E.~V.}\binits{E.~V.}}
(\byear{2013}).
\bhowpublished{Supplement to ``Note on distribution free testing for
  discrete distributions.'' DOI:\doiurl{10.1214/13-AOS1176SUPP}.}
\bptok{imsref}%
\end{bmisc}
\endbibitem


%b18 #&#
\bibitem[\protect\citeauthoryear{Khmaladze}{1993}]{Khm93}
\begin{barticle}[mr]
\bauthor{\bsnm{Khmaladze},~\bfnm{{\`E}.~V.}\binits{{\`E}.~V.}}
(\byear{1993}).
\btitle{Goodness of fit problem and scanning innovation martingales}.
\bjournal{Ann. Statist.}
\bvolume{21}
\bpages{798--829}.
\bid{doi={10.1214/aos/1176349152}, issn={0090-5364}, mr={1232520}}
\bptok{imsref}%
\end{barticle}
\endbibitem

%b19 #&#
\bibitem[\protect\citeauthoryear{Kolmogorov}{1933}]{Kol}
\begin{bmisc}[auto:STB|2013/10/14|10:36:11]
\bauthor{\bsnm{Kolmogorov},~\bfnm{A.~N.}\binits{A.~N.}}
(\byear{1933}).
\bhowpublished{Sulla determinazione empirica di una legge di
  distribuzione, Giornale dell'Istituto Italiano degli Attuari; see also
\textsc{Kolmogorov,~A.~N.} (1992). \textit{Selected Works. Vol.~II}: \textit{Probability theory and mathematical statistics}.
Kluwer Academic, Dordrecht}.
\bid{mr={1153022}}
\bptok{imsref}%
\end{bmisc}
\endbibitem

%%%\bibitem[\protect\citeauthoryear{Kolmogorov}{1992}]{K92}
%%%\begin{bbook}[mr]
%%%\bauthor{\bsnm{Kolmogorov},~\bfnm{A.~N.}\binits{A.~N.}}
%%%(\byear{1992}).
%%%\btitle{Selected Works. {V}ol. {II}: Probability theory and mathematical statistics}.
%%%\bpublisher{Kluwer Academic}, \blocation{Dordrecht}.
%%%\bid{mr={1153022}}
%%%\bptok{imsref}%
%%%\end{bbook}
%%%\endbibitem

%b20 #&#
\bibitem[\protect\citeauthoryear{Owen}{2001}]{Owen}
\begin{bbook}[auto:STB|2013/10/14|10:36:11]
\bauthor{\bsnm{Owen},~\bfnm{A.}\binits{A.}}
(\byear{2001}).
\btitle{Empirical Likelihood}.
\bpublisher{Chapman \& Hall}, \blocation{Boca Raton, FL}.
\bptok{imsref}%
\end{bbook}
\endbibitem

%b21 #&#
\bibitem[\protect\citeauthoryear{Pearson}{1900}]{Pear}
\begin{barticle}[auto:STB|2013/10/14|10:36:11]
\bauthor{\bsnm{Pearson},~\bfnm{K.}\binits{K.}}
(\byear{1900}).
\btitle{On the criterion that a given system of deviations from the probable in
  the case of a correlated system of variables is such that it can be
  reasonably supposed to have arisen from random sampling}.
\bjournal{Philosophical Magazine}
\bvolume{50}
\bpages{157--175}.
\bnote{Reprinted in Karl Pearson's Early Statistical Papers, 1948, 339--357,
Cambridge Univ. Press, Cambridge.}
\bptok{imsref}%
\end{barticle}
\endbibitem

%b22 #&#
\bibitem[\protect\citeauthoryear{Rao}{1965}]{Rao}
\begin{bbook}[mr]
\bauthor{\bsnm{Rao},~\bfnm{C.~Radhakrishna}\binits{C.~R.}}
(\byear{1965}).
\btitle{Linear Statistical Inference and Its Applications}.
\bpublisher{Wiley}, \blocation{New York}.
\bid{mr={0221616}}
\bptnote{check year}%
\bptok{imsref}%
\end{bbook}
\endbibitem

%b23 #&#
\bibitem[\protect\citeauthoryear{Rosenblatt}{1952}]{Ros}
\begin{barticle}[mr]
\bauthor{\bsnm{Rosenblatt},~\bfnm{Murray}\binits{M.}}
(\byear{1952}).
\btitle{Remarks on a multivariate transformation}.
\bjournal{Ann. Math. Statistics}
\bvolume{23}
\bpages{470--472}.
\bid{issn={0003-4851}, mr={0049525}}
\bptok{imsref}%
\end{barticle}
\endbibitem

%%b24 #&#
%(\byear{1997}).

\bibitem[\protect\citeauthoryear{Smirnov}{1937}]{Sil}
\begin{barticle}[auto:STB|2013/12/09|07:59:19]
\bauthor{\bsnm{Smirnov},~\bfnm{N.~V.}\binits{N.~V.}}
(\byear{1937}).
\btitle{On the distribution of $\omega^2$-test of Mises}.
\bjournal{Mat. Sb.}
\bvolume{2}
\bpages{973--993}
\bnote{(in Russian)}.
\bptok{imsref}%
\end{barticle}
\endbibitem


%b25 #&#
\bibitem[\protect\citeauthoryear{Stigler}{1999}]{Stig}
\begin{bbook}[mr]
\bauthor{\bsnm{Stigler},~\bfnm{Stephen~M.}\binits{S.~M.}}
(\byear{1999}).
\btitle{Statistics on the Table}.
\bpublisher{Harvard Univ. Press}, \blocation{Cambridge, MA}.
\bid{mr={1712969}}
\bptok{imsref}%
\end{bbook}
\endbibitem

%b26 #&#
\bibitem[\protect\citeauthoryear{van~der Vaart}{1998}]{vdV}
\begin{bbook}[mr]
\bauthor{\bparticle{van~der} \bsnm{Vaart},~\bfnm{A.~W.}\binits{A.~W.}}
(\byear{1998}).
\btitle{Asymptotic Statistics}.
\bpublisher{Cambridge Univ. Press}, \blocation{Cambridge}.
\bid{mr={1652247}}
\bptok{imsref}%
\end{bbook}
\endbibitem

%b27 #&#
\bibitem[\protect\citeauthoryear{Wald and Wolfowitz}{1939}]{WalWol}
\begin{barticle}[auto:STB|2013/10/14|10:36:11]
\bauthor{\bsnm{Wald},~\bfnm{A.}\binits{A.}} \AND
  \bauthor{\bsnm{Wolfowitz},~\bfnm{J.}\binits{J.}}
(\byear{1939}).
\btitle{Confidence limits for continuous distribution function}.
\bjournal{Ann. Math. Statistics}
\bvolume{10}
\bpages{105--118}.
\bptok{imsref}%
\end{barticle}
\endbibitem

\end{thebibliography}
\end{document}